\documentclass[a4paper,12pt]{article}

\title{
\LARGE{A proof of simultaneous linearization\\
 with a polylog estimate
}
}
\author{
Tomoki Kawahira
\\ \small{Nagoya University}
}

\usepackage{enumerate}
\usepackage{graphicx}
\usepackage{amsmath}
\usepackage{amscd}
\usepackage{amssymb}
\usepackage{amsfonts}
\usepackage{theorem}
\newtheorem{thm}{Theorem}[section]
\newtheorem{prop}[thm]{Proposition}
\newtheorem{lem}[thm]{Lemma}

\theorembodyfont{\rmfamily}

\newtheorem{pf}{Proof.}

\newcommand{\thmref}[1]{Theorem \ref{#1}}
\newcommand{\propref}[1]{Proposition \ref{#1}}

\newcommand{\lemref}[1]{Lemma \ref{#1}}

\newcommand{\B}{\mathbb{B}}
\newcommand{\C}{\mathbb{C}}

\newcommand{\Chat}{\hat{\mathbb{C}}}
\newcommand{\Cbar}{\bar{\mathbb{C}}}

\newcommand{\D}{\mathbb{D}}

\newcommand{\V}{\mathbb{V}}
\newcommand{\N}{\mathbb{N}}

\newcommand{\abs}[1]{{\left| #1 \right|}}
\newcommand{\kakko}[1]{{\left( #1 \right)}}
\newcommand{\skakko}[1]{{\left\{ #1 \right\}}}

\newcommand{\QED}{\hfill $\blacksquare$}
\newcommand{\ee}{~=~}
\newcommand{\dee}{~:=~}
\newcommand{\gee}{~\ge~}
\newcommand{\lee}{~\le~}

\newcommand{\Li}{\mathrm{Li}}

\newcommand{\rp}{\mathrm{Re}\,}

\newcommand{\al}{\alpha}

\newcommand{\lam}{\lambda}

\newcommand{\s}{\sigma}
\newcommand{\e}{\epsilon}
\newcommand{\cc}{\circ}
\newcommand{\fe}{f_\epsilon}
\newcommand{\gep}{g_\epsilon}

\newcommand{\Kfc}{K_f^\circ}
\newcommand{\Kgc}{K_g^\circ}

\newcommand{\taue}{\tau_\epsilon}
\newcommand{\elle}{\ell_\epsilon}

\newcommand{\Lame}{\varLambda_\epsilon}

\newcommand{\Be}{B_\epsilon}
\newcommand{\ue}{u_\epsilon}
\newcommand{\be}{b_\epsilon}
\newcommand{\Ve}{\V_\epsilon}
\newcommand{\Ne}{N_\epsilon}
\newcommand{\Phie}{\Phi_\epsilon}
\newcommand{\Psie}{\Psi_\epsilon}
\newcommand{\Phietilde}{\tilde{\Phi}_\epsilon}

\begin{document}

\maketitle

\begin{abstract}
We give an alternative proof of simultaneous linearization recently shown by T.Ueda, which connects the Schr\"oder equation and the Abel equation analytically. Indeed, we generalize Ueda's original result so that we may apply it to the parabolic fixed points with multiple petals. As an application, we show a continuity result on linearizing coordinates in complex dynamics. 
\end{abstract}

\section{Introduction}\label{sec_01}
Let us start with a work out example to explain a motivation to consider the simultaneous linearization theorem.

\paragraph{Cauliflowers.} 
In the family of quadratic maps, the simplest parabolic fixed point is given by $g(w)=w+w^2$. (Whose Julia set is called the \textit{cauliflower}.) Now we consider its perturbation of the form $f(w)=\lam w +w^2$ with $\lam \nearrow 1$. According to \cite[\S 8 and \S 10]{MiBook}, we have the following fact:

\begin{prop}[K\"onigs and Fatou coordinates]\label{prop_Kf_Kg}
Let $K_f$ and $K_g$ be the filled Julia sets of $f$ and $g$. Then we have the following:
\begin{enumerate}
\item There exists a unique holomorphic branched covering map $\phi_f:\Kfc \to \C$ satisfying the Schr\"oder equation $\phi_f(f(w))=\lam \phi_f(w)$ and $\phi_f(0)=\phi_f(-\lam/2)-1=0$. $\phi_f$ is univalent near $w=0$.
\item There exists a unique holomorphic branched covering map $\phi_g:\Kgc \to \C$ satisfying the Abel equation $\phi_g(g(w))=\phi_g(w)+1$ and $\phi_g(-1/2)=0$. $\phi_g$ is univalent on a disk $|w+r|<r$ with small $r>0$.
\end{enumerate}
\end{prop}
Note that $-\lam/2$ and $-1/2$ are the critical points of $f$ and $g$ respectively.

\paragraph{Observation.}
Set $\tilde{w}=\phi_f(w)$. Now the proposition above claims that the action of $f|_{\Kfc}$ is semiconjugated to $\tilde{w} \mapsto \lam \tilde{w}$ by $\phi_f$. Let us consider a M\"obius map $W=S_f(\tilde{w})=\lam(\tilde{w}-1)/(\lam-1)\tilde{w}$ that sends $\skakko{0, 1, \lam}$ to $\skakko{\infty,0,1}$ respectively. By taking conjugation by $S_f$, the action of $\tilde{w} \mapsto \lam \tilde{w}$ is viewed as $W \mapsto W/\lam+1$. Let us set $W=\Phi_f(w):=S_f \cc \phi_f(w)$. Now we have
$$
\Phi_f(f(w)) \ee \Phi_f(w)/\lam + 1 \text{~~~and~~~} \Phi_f(-\lam/2) \ee 0.
$$
in total. On the other hand, by setting $W=\Phi_g(w):=\phi_g(w)$, we can see the action of $g|_{\Kgc}$ as $W \mapsto W+1$. Thus we have
$$
\Phi_g(g(w)) \ee \Phi_g(w) + 1 \text{~~~and~~~} \Phi_g(-1/2) \ee 0.
$$
If $\lam$ tends to $1$, that is, $f \to g$, the semiconjugated action in $W$-coordinate converges uniformly on compact sets. Now it would be natural if $\Phi_f$ tends to $\Phi_g$. However, as one can see by referring the proof of the proposition in \cite[\S 8 and \S 10]{MiBook}, $\phi_f$ and $\phi_g$ are given in completely different ways thus we cannot conclude the convergence of $\Phi_f \to \Phi_g$ a priori.

\begin{figure}[htbp]
\label{fig_Phi}
\centering{\vspace{4cm}
}
\caption{Semiconjugation inside the filled Julia sets of Cauliflowers.}
\end{figure}

But there is another evident that supports this observation. Figure 1 shows the equipotential curves of $\phi_f$ and $\phi_g$ in the filled Julia sets. We can find similar patterns and it seems one converges to another. Actually, we have the following:

\begin{thm}\label{thm_convergence}
For any compact set $E \subset \Kgc$, 
\begin{enumerate}[(1)]
\item $E \subset \Kfc$ for all $f \approx g$; and
\item $\Phi_f \to \Phi_g$ uniformly on $E$ as $f \to g$.
\end{enumerate}
\end{thm}
Here $f \approx g$ means that $f$ is sufficiently close to $g$, equivalently, $\lam$ sufficiently close to $1$. (See \cite[Theorem 5.5]{Ka} for more generalized version of this proposition.) The proof of this theorem is given in Section 5, by using the \textit{simultaneous linearization theorem}.

\section{Simultaneous linearization}\label{sec_02}

In this section we state the simultaneous linearization theorem. We first generalize the cauliflower setting above:

\paragraph{Perturbation of parabolics.}
Let $f$ be an analytic map defined on a neighborhood of $0$ in $\Cbar$ which is tangent to identity at $0$. That is, $f$ near $0$ is of the form
$$
f(w) \ee w + A w^{m+1}+O(w^{m+2})
$$
where $A \neq 0$ and $m \in \N$. By taking a linear coordinate change $w \mapsto A^{1/m}w$, we may assume that $A=1$. In the theory of complex dynamics such a germ appears when we consider iteration of local dynamics near the parabolic periodic points, and plays very important roles. (See \cite[\S 10]{MiBook} and \cite{Sh} for example.) Now we consider a perturbation $\fe \to f$ of the form
$$
\fe(w) \ee \Lame w \, (1+ w^{m}+O(w^{m+1}))
$$
with $\Lame \to 1$ as $\e \to 0$. By taking branched coordinate changes $z=-\Lame^{m}/(m w^m)$ and setting $\taue:=\Lame^{-m}$, we have 
\begin{align*}
\fe(z) & \ee \taue  z +1 +O(|z|^{-1/m}) \\
\longrightarrow 
f_0(z) & \ee z +1 +O(|z|^{-1/m})
\end{align*}
uniformly near $w=\infty$ on the Riemann sphere $\Chat$. The simultaneous linearization theorem will give partially linearizing coordinates of $\fe$ that depend continuously on $\e$ when $\taue \to 1$ non-tangentially to the unit circle.

Let us formalize non-tangential accesses to $1$ in the complex plane: After C.McMullen, for a variable $\tau \in \C$ converging to $1$, we say $\tau \to 1$ \textit{radially} (or more precisely, \textit{$\al$-radially}) if $\tau$ satisfies $|\arg(\tau-1)| \le \alpha$ for some fixed $\al \in [0, \pi/2)$. 

\paragraph{Ueda's modulus.} 
Let us consider a continuous family of complex numbers $\skakko{\taue}$ with $\e \in [0,1]$ such that $|\taue|\ge 1$ and $\taue ~\to~ 1$ $\al$-radially as $\e \to 0$. For simplicity we assume that $\taue = 1$ iff $\e=0$. Set $\elle(z):=\taue z +1$, which is an isomorphism of the Riemann sphere $\Chat$. If $\e>0$, then $\be:=1/(1-\taue)$ is the repelling fixed point of $\elle$ with $\elle(z)-\be=\taue(z-\be)$. Thus the function
$$
N_\e(z) \dee |z-\be|-|\be|
$$
satisfies the uniformly increasing property 
$$
N(\elle(z)) \ee |\taue| N(z)+ \frac{|\taue|-1}{|\taue-1|} \gee |\taue| N(z) + \cos \al.
$$
Similarly, if $\e=0$, the function 
$$
N_0(z)\dee \sup\skakko{\rp (e^{i\theta}z)~:~|\theta|<\alpha}
$$
also has the corresponding property
$$
N_0(\ell_0(z)) \gee N_0(z) + \cos \al.
$$
In both cases, set 
$$
\Ve(R) \dee \skakko{z \in \C: N_\e(z) \ge R}
$$
for $R > 0$. One can check that $\Ne(z) \le |z|$ and $\Ve(R) \subset \B(R):=\skakko{z \in \C:|z| \ge R}$ for all $\e \in [0,1]$. 

We will establish:

\begin{thm}[Simultaneous linearization]\label{thm_ueda}
Let $\skakko{\fe: \e \in [0,1] }$ be a family of holomorphic maps on $\B(R)$ such that as $\e \to 0$ we have the uniform convergence on compact sets of the form
\begin{align*}
\fe(z) & \ee \taue z +1 +O(1/|z|^\sigma) \\
\longrightarrow 
f_0(z) & \ee z +1 +O(1/|z|^\sigma)
\end{align*}
for some $\sigma \in (0,1]$ and $\taue \to 1$ $\al$-radially. If $R \gg 0$, then: 
\begin{enumerate}[(1)]
\item For any $\e \in [0,1]$ there exists a holomorphic map $\ue:\V_\e(R) \to \C$ such that
$$
\ue(\fe(z)) \ee \taue \ue(z)+1.
$$
\item For any compact set $K$ contained in $\Ve(R)$ for all $\e \in [0,1]$, $\ue \to u_0$ uniformly on $K$.
\end{enumerate}
\end{thm}

This theorem is a mild generalization of Ueda's theorem in \cite{Ue} that deals with the case of $\sigma=1$. (See also \cite{Ue2}.) 
This plays a crucial role to show the continuity of tessellation of the filled Julia set for hyperbolic and parabolic quadratic maps. See \cite{Ka}. C.McMullen showed that there exist quasiconformal linearizations with much wider domain of definition. Indeed, $\taue \to 1$ may be tangent to the unit circle (\textit{horocyclic} in his term). See \cite[\S 8]{Mc}.

\paragraph{Remark on the domain of convergence.}
We can take such a compact subset $K$ above in  
\begin{align*}
\Pi(R) &\dee \C-\skakko{e^{\theta i}z : \rp z < R,~|\theta|\le \al} \\ 
&~\ee \skakko{z \in \C: \rp(z-R') \ge |z-R'| \sin\alpha},
\end{align*}
which is a closed sector at $z=R'=R/(\cos \alpha)>0$. In fact, for any $R>0$ and $\e \in [0,1]$, $\Pi(R)$ is contained in $\Ve(R)$. One can check it as follows: Now the complement of $\Ve(R)$ is contained in $\skakko{e^{\theta i}z : \rp z < R,~\theta=\arg(-\be)}$. Since $|\arg(-\be)| \le \alpha$, we have the claim.

In the next section we give a proof of this theorem that is also an alternative proof of Ueda's simultaneous linearization when $\sigma=1$. His original proof given in \cite{Ue} uses a technical difference equation but it makes the proof beautiful and the statement a little more detailed. Here we present a simplified proof based on the argument of \cite[Lemma 10.10]{MiBook} (its idea can be traced back at least to Leau's work on the Abel equation \cite{L}) and an estimate on polylogarithm functions given in Section 4. 

\section{Proof of the theorem}\label{sec_03}
Let us start with a couple of lemmas. Set $\delta := (\cos \al)/2>0$. We first check:
\begin{lem}\label{lem_ueda1}
If $R \gg 0$, there exists $M>0$ such that $|\fe(z)-(\taue z+1)| \le M/|z|^\sigma$ on $\B(R)$ and $N_\e(\fe(z)) \ge N_\e(z)+ \delta$ on $\V_\e(R)$ for any $\e \in [0,1]$.
\end{lem}  
\begin{pf}
The first inequality and the existence of $M$ is obvious. By replacing $R$ by a larger one, we have $|\fe(z)-(\taue z+1)| \le M/R^\sigma < \delta$ on $\B(R)$. Then
$$
\Ne(\fe(z)) \gee \Ne (\elle(z)) -\delta \gee \Ne(z)+\delta.
$$   
\QED
\end{pf}
Let us fix such an $R\gg 0$. Then the lemma above implies that $\fe(\Ve(R)) \subset \Ve(R)$. Moreover, since $\Ne(z) \le |z|$, we have
\begin{equation}\label{eq_1}
|\fe^n(z)| \gee \Ne(\fe^n(z)) \gee \Ne(z)+n\delta \gee R+n\delta 
~\to~ \infty. \tag{2.1}
\end{equation}
Thus $\Ve(R)$ is contained in the basin of infinity and uniformly attracted to $\infty$ in spherical metric of $\Chat$. In particular, this convergence to $\infty$ is uniform on $\Pi(R)$ for any $\e$.

Next we show a key lemma for the theorem:
\begin{lem}\label{lem_key}
There exists $C >0$ such that for any $\e \in [0,1]$ and $z_1,~z_2 \in \B(2S)$ with $S > R$, we have:
$$
\abs{\frac{\fe(z_2)-\fe(z_1)}{z_2-z_1}-\taue} \lee \frac{C}{S^{1+\sigma}}.
$$
\end{lem}

\paragraph{Proof.}
Set $\gep(z):=\fe(z)-(\taue z+1)$. For any $z \in \B(2S)$ and $w \in \D(z,S):=\skakko{w:|w-z|<S}$, we have $|w|>S$. This implies $|\gep(w)| \le M/|w|^\sigma < M/S^\sigma$ and thus $\gep$ maps $\D(z,S)$ into $\D(0,M/S^\sigma)$. By the Cauchy integral formula (or the Schwarz lemma), it follows that $|\gep'(z)| \le (M/S^\sigma)/S=M/S^{1+\sigma}$ on $\B(2S)$. 

Let $[z_1,z_2]$ denote the oriented line segment from $z_1$ to $z_2$. If $[z_1,z_2]$ is contained in $\B(2S)$, the inequality easily follows by
$$
|\gep(z_2)-\gep(z_1)| \ee \abs{\int_{[z_1,z_2]} \gep'(z) dz } 
 \lee \int_{[z_1,z_2]} |\gep'(z)| |dz| 
 \lee \frac{M}{S^{1+\sigma}}|z_2-z_1|
$$
with $C:=M$. Otherwise we have to take a roundabout way to get the estimate. Let us consider a circle with diameter $[z_1,z_2]$. Then $[z_1,z_2]$ cut the circle into two semicircles, and at least one of them is contained in $\B(2S)$. Let $\skakko{z_1,z_2}$ denote the one. Then 
$$
|\gep(z_2)-\gep(z_1)| \ee \abs{\int_{\skakko{z_1,z_2}} \gep'(z) dz } 
 \lee \int_{\skakko{z_1,z_2}} |\gep'(z)| |dz| 
 \lee \frac{M}{S^{1+\sigma}}\cdot \frac{\pi}{2}|z_2-z_1|
$$
and the lemma holds by setting $C:=M \pi/2~(>M)$ for any $z_1, ~z_2 \in \B(2S)$.
\QED

\paragraph{Proof of \thmref{thm_ueda}.}
Set $z_n:=\fe^n(z)$ for $z \in \Ve(2R)$. Note that such $z_n$ satisfies $|z_n| \ge \Ne(z_n) \ge 2R+n\delta$ by (\ref{eq_1}). Now we fix $a \in \Ve(2R)$ and define $\phi_{n,\e}=\phi_n:\Ve(2R) \to \C~(n \ge 0)$ by
$$
\phi_n(z) \dee \frac{z_n -a_n}{\taue^n}.
$$
For example, one can take such an $a$ in $\Pi(2R)$ independently of $\e$. Then we have
$$
\abs{\frac{\phi_{n+1}(z)}{\phi_n(z)}-1}
 \ee \abs{ \frac{z_{n+1} - a_{n+1}}{\taue(z_{n} - a_{n})} -1} \\
 \ee \frac{1}{|\taue|} 
 \cdot \abs{ \frac{\fe(z_{n}) - \fe(a_{n})}{z_{n} - a_{n}} -\taue}. 
$$
We apply \lemref{lem_key} with $2S=2R+n\delta$. Since $z_n,~a_n \in \Ve(2S) \subset \B(2S)$, we have 
$$
\abs{\frac{\phi_{n+1}(z)}{\phi_n(z)}-1} \lee \frac{C}{|\taue|(R+n\delta/2)^{1+\sigma}} \lee \frac{C'}{(n+1)^{1+\sigma}},
$$
where $C'=2^{1+\sigma}C/\delta^{1+\sigma}$ and we may assume $R>\delta/2$. Set $P:=\prod_{n \ge 1}(1+C'/n^{1+\sigma})$. Since $|\phi_{n+1}(z)/\phi_n(z)|\le 1+C'/(n+1)^{1+\sigma}$, we have
$$
|\phi_n(z)| \ee 
\abs{
\frac{\phi_{n}(z)}{\phi_{n-1}(z)}
} \cdots \abs{
\frac{\phi_{1}(z)}{\phi_0(z)}
} \cdot |\phi_0(z)| \lee P|z-a|.
$$
Hence 
$$
|\phi_{n+1}(z)-\phi_n(z)|  \ee 
\abs{\frac{\phi_{n+1}(z)}{\phi_n(z)}-1} \cdot |\phi_n(z)|
\lee \frac{C'P}{(n+1)^{1+\sigma}}\cdot |z-a|.
$$
This implies that $\phi_\e=\phi_0+(\phi_1-\phi_0)+\cdots=\lim \phi_n$ converges uniformly on compact subsets of $\Ve(2R)$ and for all $\e \in [0,1]$. The univalence of $\phi_\e$ is shown in the same way as \cite[Lemma 10.10]{MiBook}. 

Next we claim that $\phi_\e(\fe(z))=\taue \phi_\e(z)+\Be$ with $\Be \to 1$ as $\e \to 0$. One can easily check that $\phi_n(\fe(z))=\taue \phi_{n+1}(z)+B_n$ where 
$$
B_n \ee \frac{a_{n+1}-a_n}{\taue^n} 
\ee 
\frac{(\taue-1)a_n}{\taue^n}+
\frac{1+\gep(a_n)}{\taue^n}.
$$
When $\taue=1$, $B_n$ tends to $1$ since 
$$
|\gep(a_n)| \lee \frac{M}{|a_n|^\sigma} \lee \frac{M}{(2R+n\delta)^\sigma}
\lee \frac{M}{(n\delta)^\sigma} ~\to~ 0.
$$ 
When $|\taue| >1$, the last term of the equation on $B_n$ above tends to 0. For $n \ge 1$, we have
$$
a_n \ee \taue^n a ~+~\frac{\taue^n-1}{\taue-1}~+~\sum_{k=0}^{n-1}\taue^{n-1-k}\gep(a_{k}).
$$
Thus
$$
\frac{(\taue-1)a_n}{\taue^n} \ee
(\taue -1) \kakko{a+ \frac{\gep(a)}{\taue}
+\sum_{k=1}^{n-1}\frac{\gep(a_{k})}{\taue^{k+1}}
}+1-\frac{1}{\taue^n}.
$$
By the inequality on $|\gep(a_n)|$ above, we have 
$$
\abs{(\taue -1) \sum_{k=1}^{n-1}\frac{\gep(a_{k})}{\taue^{k+1}} }
\lee \frac{M}{\delta^\sigma} \frac{|\taue-1|}{|\taue|}\sum_{k=1}^{n-1} \frac{1}{k^\sigma |\taue|^k} 
\lee \frac{M}{2\delta^{1+\sigma}} (1-\frac{1}{|\taue|})\Li_\sigma(\frac{1}{|\taue|})
$$
where we used the inequality 
$$
|\taue -1| \lee \frac{\rp \taue -1}{\cos \al} \lee \frac{|\taue|-1}{2\delta}
$$
that comes from the radial convergence. By \propref{prop_polylog} in the next section, $B_n$ converges to some $B_\e$. More precisely, if we set $|\taue|=e^{L}$, then $\taue-1=O(L)$ and one can check that $B_\e=1+O(L^{\sigma/(1+\sigma)})$.

Finally, $\ue(z):=\phi_\e(z)/\Be$ gives a desired holomorphic map (with $R$ in the statement replaced by $2R$). \QED

\paragraph{Remarks.}
\begin{itemize}
\item 
When $\sigma=1$, we have 
$$
\abs{\sum_{k=1}^{n-1}\frac{\gep(a_{k})}{\taue^{k+1}} }
\lee \frac{M}{\delta|\taue|} \sum_{k=1}^{n-1} \frac{1}{k|\taue|^k} 
\lee -\frac{M}{\delta} \log (1-\frac{1}{|\taue|})
$$
and this implies that $B_\e=1+O(L|\log L|)$ if we set $|\taue|=e^{L}$. This fact is consistent with the result in \cite{Ue}.
\item By this proof, if $\skakko{\fe(z)}$ analytically depends on $\e$, then $\skakko{\Be}$ and $\skakko{\ue(z)}$ do the same for fixed $a$ in $\Pi(2R)$. 
\item It is not difficult to check that $u_\e(z)=z(\Be^{-1}+o(1))$ as $z \to \infty$ within $\Ve(R)$. (It is well-known that if $f_0(z)=z+1+a_0/z+\cdots$ then the Fatou coordinate is of the form $u_0(z)=z-a_0 \log z+O(1)$. See \cite{Sh}.)
\end{itemize}

\section{An estimate on polylogarithm functions}\label{sec_04}
We define the \textit{polylogarithm function} of exponent $s \in \C$ by
$$
\Li_s(z) \dee \sum_{n=1}^\infty \frac{z^n}{n^s}.
$$
This function makes sense when $|z| < 1$ and $\s:=\rp s >0$ and it is a holomorphic function of $z$. In particular, if $\rp s >1$ the function tends to $\zeta(s)$ as $z \to 1$ within the unit disk. In the following we consider the behavior of $\Li_s(z)$ as $z \to 1$ within the unit disk when $0<\s \le 1$. We claim:

\begin{prop}\label{prop_polylog}
Suppose $0<\rp s=\s \le 1$ and $z \to 1$ with $|z|<1$. Set $\e:=1-|z| $. Then there exists a uniform constant $C$ independent of $s$ such that 
$$
|\Li_s(z)| \lee C \e^{-\frac{1}{1+\sigma}}
$$
as $z \to 1$. In particular, we have
$$
|(z-1) \, \Li_s(z)| \lee C \e^{\frac{\sigma}{1+\sigma}} ~\to~ 0
$$
if $z \to 1-0$ along the real axis.
\end{prop}

\newcommand{\nsum}{\sum_{n=1}^\infty}
\newcommand{\ksum}{\sum_{k=1}^n}

\paragraph{Proof.}
Clearly 
$
|\Li_s(z)| \lee \sum_{n=1}^\infty |z|^n/n^\s
$
so it is enough to consider the sum
$$
S \dee \nsum \frac{1}{n^\s}\cdot \lam^n 
$$
where $\lam:=|z|=1-\e$. Let $S_n$ be the partial sum to the $n$th term. By the H\"older inequality, we have 
$$
S_n \lee 
\kakko{
\ksum \frac{1}{k^{\s p}}
}^{\frac{1}{p}}
\kakko{
\ksum \lam^{kq}
}^{\frac{1}{q}} 
$$
for any $p, q > 1$ with $1/p+1/q=1$. Now let us set $p:=1/\s+1 \ge 2$ (then $1<q=1+\s \le 2$). Since $\s p=1+\s>1$, the first sum is uniformly bounded as follows:
$$
\ksum \frac{1}{k^{\s p}} \lee 1+\int_1^\infty \frac{1}{x^{1+\s}} dx
\ee 1+\frac{1}{\s} = p. 
$$
On the other hand, for the second sum, we still have $0<\lam^q<1$ and thus 
$$
\ksum \lam^{kq} \lee \frac{\lam^q}{1-\lam^q} 
\ee \frac{1}{q \e}(1+o(1)) \lee \frac{2}{q \e}
$$
when $\e \ll 1$. Hence we have the following uniform bound:
$$
S_n \lee  p^{\frac{1}{p}} 
\kakko{\frac{2}{q \e}}^{\frac{1}{q}}
\lee 2 \kakko{\frac{p^{\frac{1}{p}}}{q^{\frac{1}{q}}}} \e^{-\frac{1}{q}}.
$$
One can easily check that $1 \le x^{\frac{1}{x}} \le e^{\frac{1}{e}}=1.44467\cdots$ for $x \ge 1$. Thus 
$$
S \lee 2 e^{\frac{1}{e}} \e^{-\frac{1}{q}} \ee  2 e^{\frac{1}{e}} \e^{-\frac{1}{1+\s}}
$$
when $\e \ll 1$ and we have the desired estimate with $C=2 e^{\frac{1}{e}}< 3$. The last inequality of the statement follows by:
$$
|(z-1) \, \Li_s(z)| \lee  C \e^{1-\frac{1}{q}} \ee C \e^{\frac{1}{p}}
 \ee C \e^{\frac{\s}{1+\s}}.
$$
(Indeed, $|(z-1) \, \Li_s(z)| = O(\e^{\frac{\s}{1+\s}})$ if $z \to 1$ radially.)
\QED

\section{Application: Proof of \thmref{thm_convergence}}\label{sec_05}

As an application of \thmref{thm_ueda}, we give a proof of \thmref{thm_convergence}. Though \thmref{thm_convergence} only deals with the simplest parabolic fixed point and its simplest perturbation, one can easily extend the result to general parabolic cycles with multiple petals and their ``non-tangential" perturbations. 

\paragraph{Proof of \thmref{thm_convergence}.}
Let us take an general expression $f_\lam (w)=\lam w+w^2$ with $0<\lam \le 1$ (thus $f_1=g$). By looking the action of $f_\lam$ through a new coordinate $z=\chi_\lam(w)=-\lam^2/w$, we have 
$$
\chi_\lam\cc f_\lam \cc \chi_\lam^{-1}(z) \ee z/\lam + 1 + O(1/z)
$$
near $\infty$. Now we can set $\taue:=1/\lam=1+\e$ and $\fe :=\chi_\lam\cc f_\lam \cc \chi_\lam^{-1}$ to have the same setting as \thmref{thm_ueda}. We consider that $f$ and $g$ are parameterized by $\lam$ or $\e$. (It is convenient to use both parameterization.) Note that $\Pi(R)=\skakko{\rp z  \ge R}$ in this case. By the same argument as \lemref{lem_ueda1}, we can check that $\rp \fe(z) \ge \rp z + 1/2$ if $z \in \Pi(R)$ and $R \gg 0$. In particular, we have $\fe(\Pi(R)) \subset \Pi(R)$ for $R \gg 0$.

Let us show (1): For any compact $E \subset \Kgc$ and small $r>0$, there exists $N \gg 0$ such that $g^N(E) \subset P_r=\skakko{|w+r| \le r}$. (For instance, one can show this fact by existence of the Fatou coordinate.) By uniform convergence, we have $f^N(E) \subset P_r$ for all $f \approx g$. To show $E \subset \Kfc$, it is enough to show that $f(P_r) \subset P_r$ for all $f \approx g$. Since $\chi_\lam(P_r)=\Pi(R)$ for some $R \gg 0$, we have $\fe(\Pi(R)) \subset \Pi(R)$ independently of $\e$. This is equivalent to $f_\lam(P_r) \subset P_r$ in a different coordinate. Thus we have (1).

Next let us check (2): Set $\Phi_\e :=\Phi_f$ and $\Phi_0:=\Phi_g$. Then we have $\Phie(f_\lam(w))=\taue \Phie(w)+1$. On the other hand, by simultaneous linearization, we have a uniform convergence $\ue \to u_0$ on $\Pi(R)$ that satisfies $\ue(\fe(z))=\taue \ue(z)+1$. By setting $\Psie(w):=\ue \cc \chi_\lam(w)$, we have $\Psie \to \Psi_0$ compact uniformly on $P_r$, and $\Psie(f_\lam(w))=\taue \Psie(w)+1$. 

We need to adjust the images of critical orbits mapped by $\Phie$ and $\Psie$. Since $g^n(-1/2) \to 0$ along the real axis, there is an $M \gg 0$ such that $g^M(-1/2)=:a_0 \in P_r$. By uniform convergence, we also have $f^M(-\lam/2)=: a_\e \in P_r$ and $a_\e \to a_0$ as $\e \to 0$. Set $b_\e:=\Psie(a_\e)$ and $c_\e:=\Phie(a_\e)$ for all $\e \ge 0$. Set also $\ell_\e(W)=\taue W +1$, then we have $c_\e=\ell_\e^M(0)=\taue^{M-1}+\cdots+\taue+1$ and $c_\e \to c_0=M$ as $\e \to 0$. When $\e >0$, we take an affine map $T_\e$ that fixes $1/(1-\taue)$ and sends $b_\e$ to $c_\e$. When $\e=0$, we take $T_0$ that is the translation by $b_0-c_0$. Then one can check that $T_\e \to T_0$ compact uniformly on the plane and $\Phietilde:=T_\e \cc \Psie$ satisfies $\Phietilde \to \tilde{\Phi}_0$ on any compact sets of $P_r$. Moreover, $\Phietilde$ still satisfies $\Phietilde(f_\lam(w))=\taue \Phietilde(w)+1$ and the images of the critical orbit by $\Phie$ and $\Phietilde$ agree. Finally by uniqueness of $\phi_f$ and $\phi_g$, one can check that $\Phie=\Phietilde$ on $P_r$. 

Since
$$
\Phi_f(w) \ee \ell_\e^{-N}\cc \Phietilde(f^N(w))
~\longrightarrow ~ \ell_0^{-N}\cc \tilde{\Phi}_0(g^N(w)) \ee \Phi_g(w)
$$ 
uniformly on $E$, we have (2). \QED

\paragraph{Acknowledgement.}
I would like to thank T.Ueda for correspondence. This research is partially supported by Inamori Foundation and JSPS.

\end{document}